\documentclass[runningheads]{llncs}
\usepackage{amsmath,graphicx}

\def\email#1{{\normalsize\tt #1}}

\mainmatter

\title{Operator Calculus Approach to Solving Analytic Systems}
\author{Philip Feinsilver\inst{1} and Ren{\'e} Schott\inst{2}}
\institute{Southern Illinois University, 
 {Carbondale, IL. 62901, U.S.A.}
\email{pfeinsil@math.siu.edu}
\and {Universit\'e Henri Poincar\'e-Nancy 1},
{BP 239, 54506 Vandoeuvre-l\`es-Nancy, France}
\email{schott@loria.fr}
}

\begin{document} 
\maketitle
\begin{abstract}

Solving analytic systems using inversion can be implemented in a variety of ways.
One method is to use  Lagrange inversion and variations. Here we present a different approach,
based on dual vector fields.

~~~~~For a function analytic in a neighborhood of the origin in the complex plane, we associate
a vector field and its dual, an operator version of Fourier transform. The construction
extends naturally to functions of several variables.

~~~~~We illustrate with various examples and present an efficient algorithm
readily implemented as a symbolic procedure in Maple while suitable as well for
numerical computations using languages such as C or Java. 
\end{abstract}


\section{Introduction}

We introduce the operator calculus necessary to present our approach to (local) inversion of
analytic functions. It is important to note that this is different from Lagrange inversion and is 
based on the flow of a vector field associated to a given function. It appears to be 
theoretically appealing as well as computationally effective.

Acting on polynomials in $x$, define the operators

\medskip
 $\displaystyle D=\frac{d}{dx}$ and $X=\hbox{multiplication by }x$. 

\medskip\noindent
They satisfy commutation relations $[D,X]=I$, where $I$, the identity
operator, commutes with both $D$ and $X$. Abstractly, the Heisenberg-Weyl algebra is the associative
algebra generated by operators $\{A,B,C\}$ satisfying $[A,B]=C$, $[A,C]=[B,C]=0$. 
The {\sl standard\/} HW algebra is the one generated by the realization $A=D$, $B=X$, $C=I$.
An {\sl Appell system\/} is a system of polynomials $\{y_n(x)\}_{n\ge0}$ that is a basis for 
a representation of the standard HW algebra with the following properties:
\begin{enumerate}  
\item $y_n$ is of degree $n$ in $x$;
\item $D\,y_n=n\,y_{n-1}$.
\end{enumerate}
In several variables, ${\textbf x}=(x_1,\ldots,x_N)$, 
with multi-indices ${\textbf n}=(n_1,\ldots,n_N)$, 
the corresponding monomials are

\medskip
${\textbf x}^{{\textbf n}}=x_1^{n_1}x_2^{n_2}\cdots x_N^{n_N}$. 

\medskip\noindent
Denote the partial derivative operators by 
$\displaystyle D_i=\frac{\partial}{\partial x_i}$ and the corresponding
multiplication operators by $X_i$. Then $[D_j,X_i]=\delta_{ij}\,I$. 
An Appell system is a system of
polynomials $\{y_{{\textbf n}}\}$ in the variables ${\textbf x}$ such that
\begin{enumerate}
\item the top degree term of $y_{{\textbf n}}$ is a constant multiple of ${\textbf x}^{{\textbf n}}$;
\item $D_i\,y_{{\textbf n}}=n_i\,y_{{\textbf n}-{\textbf e}_i}$, 
where ${\textbf e}_i$ has all components zero except for 1 in the $i^{\rm th}$ position.
\end{enumerate}

G.-C. Rota \cite{RKO} is well-known for his {\sl umbral calculus\/} development of special polynomial sequences, called {\sl basic sequences.} 
From our perspective, these are ``canonical polynomial systems" in the sense
that they provide polynomial representations of the 
Heisenberg-Weyl algebra, in realizations different from the
standard one. Our idea  \cite{FSdvf,FSalgop} is to illustrate explicitly the r\^ole of vector fields
and their duals, using operator calculus methods for working with the latter
(in our volumes --- this viewpoint is prefigured in \cite{RKO}). 

The main feature of our approach is that the action of the vector field may be readily
calculated while the action of the dual vector field on exponentials is identical to that of 
the vector field. Then we note that acting iteratively with a vector field on 
polynomials involves the complexity
of the coefficients, while acting iteratively with the dual vector field always produces polynomials
from polynomials. So we can switch to the dual vector field for calculations.

Specifically, fix a neighborhood of 0 in ${\textbf C}$. Take an analytic function $V(z)$ defined there,
normalized to $V(0)=0$, $V'(0)=1$. Denote $W(z)=1/V'(z)$ and $U(v)$ the inverse function, i.e.,
$V(U(v))=v$, $U(V(z))=z$. Then $V(D)$ is defined by power series as an operator on polynomials
in $x$ and $[V(D),X]=V'(D)$ so that $[V(D),XW(D)]=I$. In other words, $V=V(D)$ and
$Y=XW(D)$ generate a representation of the HW algebra on polynomials in $x$. The basis for the
representation is $y_n(x)=Y^n1$, i.e., $Y$ is a {\sl raising operator\/}. And $Vy_n=n\,y_{n-1}$ so that
$V$ is the corresponding {\sl lowering operator.} The $\{y_n\}_{n\ge0}$ form a system of
{\sl canonical polynomials\/} or generalized Appell system. The operator of multiplication
by $x$ is given by $X=YV'(D)=YU'(V)^{-1}$, which is a {\sl recursion operator\/} for the system. 

We identify vector fields with first-order partial differential operators.
Consider a variable $A$ with corresponding partial differential operator $\partial_A$.
Given $V$ as above, let $\tilde Y$ be the vector field $\tilde Y=W(A)\partial_A$.
Then we observe the following identities
$$ \tilde Y\,e^{Ax}=xW(A)\,e^{Ax}=xW(D)\,e^{Ax}$$
as any operator function of $D$ acts as a multiplication operator on $e^{Ax}$.
The important property of these equalities is that $Y$ and $\tilde Y$ commute, as they
involve independent variables. So we may iterate to get

\begin{equation}\label{eq:iterate}
\exp(t\tilde Y)e^{Ax}=\exp(tY)e^{Ax}.
\end{equation}

On the other hand, we can solve for the left-hand side of this equation using the method
of characteristics. Namely, if we solve
\begin{equation}\label{eq:vf}
 \dot A=W(A)
\end{equation}
with initial condition $A(0)=A$, then for any smooth function $f$,
$$ e^{t\tilde Y}f(A)=f(A(t)).$$
Thus
$$\exp(tY)e^{Ax}=e^{xA(t)}.$$
To solve equation (\ref{eq:vf}), multiply both sides by $V'(A)$ and observe that we get
$$V'(A)\,\dot A=\frac{d}{dt}\,V(A(t))=1.$$
Integrating yields
$$V(A(t))=t+V(A)\qquad\hbox{or}\qquad A(t)=U(t+V(A)).$$
Or, writing $v$ for $t$, we have
\begin{equation}\label{eq:evf}
\exp(vY)e^{Ax}=e^{xU(v+V(A))}.
\end{equation}
We can set $A=0$ to get
$$\exp(vY)1=e^{xU(v)}$$
on the one hand while
$$ e^{vY}1=\sum_{n=0}^\infty \frac{v^n}{n!}\,y_n(x).$$
In summary, we have the expansion of the exponential of the inverse function
$$
e^{xU(v)}=\sum_{n=0}^\infty \frac{v^n}{n!}\,y_n(x)
$$
or
\begin{equation}\label{eq:expansion}
\sum_{m=0}^\infty \frac{x^m}{m!}\,(U(v))^m=\sum_{n=0}^\infty \frac{v^n}{n!}\,y_n(x).
\end{equation}
This yields an alternative approach to inversion of the function $V(z)$ rather than using Lagrange's formula. We see that the coefficient of $x^m/m!$ yields the expansion of $(U(v))^m$.
In particular, $U(v)$ itself is given by the coefficient of $x$ on the right-hand side.

Specifically, we have:
\begin{theorem}
The coefficient of $x^m/m!$ in $Y^n1$ is equal to ${\tilde Y}^nA^m\bigr|_{A=0}$, each 
giving the coefficient of $v^n/n!$ in the expansion of $U(v)^m$.
\end{theorem}
\begin{proof}
Expand both sides of  equation (\ref{eq:iterate}), using $v$ for $t$, in powers of $x$ and $v$,
and let $A=0$:
$$\sum_{n=0}^\infty \frac{v^n}{n!}\,{\tilde Y}^n\sum_{m=0}^\infty \frac{x^m}{m!}\,A^m\bigr|_{A=0}=\sum_{n=0}^\infty \frac{v^n}{n!}\,Y^n1 
$$
and compare with equation (\ref{eq:expansion}).
\end{proof}

The same idea works in several variables.

We have ${\textbf V}({\textbf z})=(V_1(z_1,\ldots,z_N),\ldots,V_N(z_1,\ldots,z_N))$ analytic
in a neighborhood of $0$ in ${\textbf C}^N$. Denote the Jacobian matrix 
$\displaystyle\biggl(\frac{\partial V_i}{\partial z_j}\biggr)$ by $V'$ and its inverse by $W$. The variables
$$ Y_i=\sum_{k=1}^N x_kW_{ki}(D)$$
commute and act as raising operators for generating the basis $y_{{\textbf n}}({\textbf x})$. Namely,
$Y_iy_{{\textbf n}}=y_{{\textbf n}+{\textbf e}_i}$. 
And $V_i=V_i({\textbf D})$, ${\textbf D}=(D_1,\ldots,D_N)$ are lowering operators: 
$V_iy_{{\textbf n}}=n_i\,y_{{\textbf n}-{\textbf e}_i}$.  

Denote $\sum_i a_ib_i$ by $a\cdot b$. With variables $A_i$ and corresponding
partials $\partial_i$, define the vector fields
$$\tilde Y_i=\sum_k W_{ki}(A)\partial_k.$$
For a vector field $\tilde Y=\sum_i W_i(A)\partial_i$, we have the identities
$$ \tilde Y\,e^{A\cdot x}=x\cdot W(A)\,e^{A\cdot x}=x\cdot W(D)\,e^{A\cdot x}.$$
The method of characteristics applies as in one variable and as in equation (\ref{eq:evf})
$$\exp(v\cdot Y)e^{A\cdot x}=e^{x\cdot U(v+V(A))}.$$
Thus, we have the expansion
\begin{equation}
\label{eq:dvfmult}
\exp\bigl(x\cdot U(v) \bigr) =  
\sum_{\textbf n}\frac{{\textbf v}^{\textbf n}}{{\textbf n}!}\,y_{{\textbf n}}({\textbf x}).
\end{equation}

In particular, the $k^{\rm th}$ component, $U_k$, of the inverse function is given 
by the coefficient of $x_k$ in the above expansion.

An important feature of our approach is that to get an expansion to a given order requires knowledge of the expansion of $W$ just to that order. 
The reason is that when iterating $xW(D)$, at step $n$ it is acting on a polynomial of
degree $n-1$, so all terms of the expansion of $W(D)$ of order $n$ or higher would
yield zero acting on $y_{n-1}$. This allows for streamlined computations. 

For polynomial systems ${\textbf V}$, 
$V'$ will have polynomial entries, and $W$ will be rational in ${\textbf z}$.
Hence raising operators will be rational functions of ${\textbf D}$, linear in ${\textbf x}$. Thus the coefficients of the expansion
of the entries $W_{ij}$ of $W$ would be computed by finite-step recurrences.

\begin{remark}
Note that to solve $V(z)=v$ for $z$ near $z_0$, with $V(z_0)=v_0$, apply the method
to $V_1(z)=V(z+z_0)-v_0$, so that $V_1(0)=0$.  The inverse is $U_1(v)=U(v+v_0)-z_0$. Then
$U(v)=z_0+U_1(v-v_0)$.
\end{remark}


\section{One-variable Case}

In this section we focus on the one-variable case.
We illustrate the method with examples, and then present 
an algorithm suitable for symbolic computation.

\begin{example} \rm    
In one variable, solving a cubic is interesting as
the expansion of $W$ can be expressed in terms of Chebyshev polynomials.

Let $V=z^3/3-\alpha z^2+z$. Then $V'=z^2-2\alpha z+1$. Thus
$$
W = \frac{1}{1-2\alpha z+z^2} = \sum_{n=0}^\infty z^n U_n(\alpha),
$$
where $U_n$ are Chebyshev polynomials of the second kind. 

Specializing $\alpha$ provides interesting cases. For example,
let $\alpha=\cos (\pi/4)$, or $V=z^3/3-z^2/\sqrt{2}+z$. Then the coefficients in the expansion
of $W$ are periodic with period 8 and, in fact,
$$ W=\frac{1+z^2+\sqrt{2}\,z}{1+z^4}.$$
The coefficient of $x$ in the polynomials $y_n$ yield the coefficients in the expansion of the inverse
$U$. Here are some polynomials starting with $y_0=1$, $y_1=x$:
\begin{eqnarray*}
y_2&=&x^2+x\,\sqrt{2},\quad y_3=x^{3} + 3\,x^{2}\,\sqrt{2} + 4\,x, \\
y_4&=&x^{4} + 6\,x^{3}\,\sqrt{2} + 22\,x^{2} + 10\,x\,\sqrt{2},\\
y_5&=&x^{5} + 10\,x^{4}\,\sqrt{2} + 70\,x^{3} + 90\,x^{2}\,\sqrt{2} + 40\,x, \\
y_6&=&x^{6} + 15\,x^{5}\,\sqrt{2} + 170\,x^{4} + 420\,x^{3}\,\sqrt{2}+ 700\,x^{2} - 140\,x\,\sqrt{2}.
\end{eqnarray*}
This gives  to order 6:
$$U(v)=\left(v+\frac{2}{3}\,v^3+\frac{1}{3}\,v^5+\cdots\right)+
\sqrt{2}\,\left(\frac{1}{2}\,v^2+\frac{5}{12}\,v^4-\frac{7}{36}\,v^6+\cdots\right).
$$
This expansion will give approximate solutions to $$z^3/3-z^2/\sqrt{2}+z-v=0$$ for $v$ near $0$.
\end{example}

\begin{example} \rm   
Inversion of the Chebyshev polynomial $T_3(z)=4z^3-3z$ 
can be used as the basis for solving general cubic equations (\cite{WIKIcubic}).

\medskip
To get started we have, with $V(z)=4z^3-3z$, 
$$W(z)=\frac{-1}{3}\,\frac{1}{1-4z^2}=\frac{-1}{3}\,
\sum_{n=0}^\infty 4^nz^{2n}.$$
So $y_1=(-1/3)x$, $y_2=(1/9)x^2$, $y_3=(-1/27)(x^3+8x)$, etc. We find
$$U(v)=-\frac{1}{3}\,v-{\frac {4}{81}}\,{v}^{3}-{\frac {16}{729}}\,{v}^{5}-{\frac {
256}{19683}}\,{v}^{7}-\cdots.$$

In this case, we can find the expansion analytically.
To solve $T_3(z)=v$, write $$T_3(\cos\theta)=\cos(3\theta)=v.$$ 
Invert to get, for integer $k$,
$\theta =(1/3)(2\pi k\pm \arccos v)$, with $\arccos$ denoting the principal branch. Then
$$z=\cos((1/3)(2\pi k\pm \arccos v)).$$ 
We want a branch with $v=0$ corresponding to $z=0$.
With $\arccos0=\pi/2$, we want the argument of the cosine to be $\pi/2+\pi l$, 
for some integer $l$. This yields the condition
$\displaystyle \frac{1}{3} = \frac{2l+1}{4k\pm1}$.
Taking $l=0$, we get $k=1$, with the minus sign. Namely, 
$$U(v)=\cos((1/3)(2\pi-\arccos v)).$$
Using hypergeometric functions (see next example) and rewriting, we find the form
$$U(v)=-\frac{1}{3}
\,\sum_{n=0}^\infty \binom{3n}{n}\,\left(\frac{4}{27}\right)^n\,
\frac{v^{2n+1}}{2n+1}.$$

If we generate the polynomials $y_n$, we can find the expansion of $U(v)^m$ to any order.

\end{example}
\begin{example}   \rm 
A similar approach is interesting for the Chebyshev polynomial $T_n(z)$.

$F(v)=\cos(\lambda(\mu\pm\arccos v))$
satisfies the hypergeometric differential equation
$$(1-v^2)\,F''-v\,F'+\lambda^2\,F=0$$
which can be written in the form
$$ [(vD_v)^2-D_v^2]F=\lambda^2\,F$$
with here $D_v$ denoting $d/dv$. For integer $\lambda$, this is the differential equation for
the corresponding Chebyshev polynomial. In general, these are {\sl Chebyshev functions\/}. As noted above,
for $F(0)=0$, we take $\mu=2\pi k$, and, as above, we require
$$\lambda=\frac{2l+1}{4k\pm1}.$$
With $F'(0)=\pm\lambda$, we have the solution
$$F(v)=\pm\lambda v\,
{}_2 F_1\left( \genfrac{}{}{0pt}{1}
{ \frac{\textstyle 1+\lambda\mathstrut}{\textstyle 2\mathstrut},
\frac{\textstyle 1-\lambda\mathstrut}{\textstyle 2\mathstrut}}
{\frac{\textstyle 3\mathstrut}{\textstyle 2\mathstrut} }
\Biggm| v^2 \right).
$$
\end{example}

\subsection{Using Maple}

For symbolic computation using Maple, one can use the Ore\_Algebra package.
\begin{enumerate}
\item First fix the degree of approximation. Expand $W$ as a polynomial to that degree.
\item Declare the Ore algebra with one variable, $x$, and one derivative, $D$.
\item Define the operator $xW(D)$ in the algebra.
\item Iterate starting with $y_0=1$ using the {\tt applyopr} command.
\item Extract the coefficient of $x^m/m!$ to get the expansion of $U(v)^m$.
\end{enumerate}

\section{Algorithm as a Matrix Computation}

Here is a matrix approach that can be implemented numerically.

Fix the order of approximation $n$. Cut off the expansion
$$W(z)=w_0+w_1z+w_2w^2+\cdots+w_kz^k+\cdots$$ 
at $w_nz^n$.

Let the matrix 
$$W=\begin{pmatrix}w_1 & w_0 &0 &\ldots&0\cr
                    w_2&w_1&w_0&\ldots&0\cr
                    \vdots&\vdots&\vdots&\ddots&\vdots\cr
                    w_{n-1}&w_{n-2}&w_{n-3}& \ldots&w_0\cr 
                    w_n&w_{n-1}&w_{n-2}&\ldots&w_1\cr\end{pmatrix}.$$

Define the auxiliary diagonal matrices 

\begin{eqnarray*}
P&=&\begin{pmatrix}1!&0&\ldots&0\cr 0&2!&\ldots&0\cr \vdots&\vdots&\ddots&\vdots\cr
                    0&0&\ldots&n!\cr\end{pmatrix}, \quad
M=\begin{pmatrix}1&0&\ldots&0\cr 0&2&\ldots&0\cr \vdots&\vdots&\ddots&\vdots\cr
                    0&0&\ldots&n\cr\end{pmatrix},\\[2pt]
Q&=&\begin{pmatrix}1/\Gamma(1)&0&\ldots&0\cr 0&1/\Gamma(2)&\ldots&0\cr \vdots&\vdots&\ddots&\vdots\cr
                    0&0&\ldots&1/\Gamma(n)\cr\end{pmatrix}.
\end{eqnarray*}
Note that $QP=M$.

Denoting $y_k(x)=\sum c_j^{(k)}x^j$, we have the recursion
$$[c_1^{(k+1)},c_2^{(k+1)},\ldots,c_n^{(k+1)}] = [c_1^{(k)},c_2^{(k)},\ldots,c_n^{(k)}] PWQ.$$
The condition $U(0)=0$ gives $y_0=1$.
Then $y_1=XW(D)y_0$ yields $y_1=w_0x$. We see that $c_0^{(k)}=0$ for $k>0$. 
We iterate as follows:

\medskip\noindent
1. Start with $w_0$ times the unit vector $[1,0,\ldots,0]$ of length $n$.

\smallskip\noindent
2. Multiply by $W$.

\smallskip\noindent
3. Iterate, multiplying on the right by $MW$ at each step.

\smallskip\noindent
4. Finally, multiply on the right by $Q$.

\medskip
The top row will give the coefficients of the expansion of $U(v)$ to order $n$.

\section{Higher-order Example}

Here is a simple $2\times2$ system for illustration.
\[V_1=z_1+z_2^2/2,\quad
V_2=z_2-z_1z_2.\]
So
$$ 
V'=\begin{pmatrix}1&z_2\cr -z_2&1-z_1\cr\end{pmatrix}\quad\hbox{and}\quad
W=\frac{1}{1-z_1+z_2^2}\,\begin{pmatrix}1-z_1&-z_2\cr z_2&1\cr\end{pmatrix}.
$$
The raising operators are
\begin{eqnarray*}
Y_1&=&\bigl(x_1(1-D_1))+x_2D_2\bigr)\,(1-D_1+D_2^2)^{-1},\\
Y_2&=&(-x_1D_2+x_2)\,(1-D_1+D_2^2)^{-1}.
\end{eqnarray*}
Expanding $\,(1-D_1+D_2^2)^{-1}=\sum\limits_{n=0}^\infty (D_1-D_2^2)^n$ yields, with $y_{00}=1$,
\begin{eqnarray*}
& y_{01}=x_2, \qquad y_{10}=x_1,\\
& y_{02}=x_2^2-x_1,  \quad y_{11}=x_2+x_1x_2, \quad  y_{20}=x_1^2.
\end{eqnarray*}
Thus
\begin{eqnarray*}
&&\exp\bigl({\textbf x}\cdot{\textbf U}({\textbf v})\bigr)=1+x_1v_1+x_2v_2 \\
&&\hspace{2.8cm}+\,(x_2+x_1x_2)v_1v_2+
(x_2^2-x_1)\,\frac{v_1^2}{2}+x_1^2\,\frac{v_2^2}{2}+\cdots,
\end{eqnarray*}
so 
\[U_1({\textbf v})= v_1-v_1^2/2+\cdots,\quad
U_2({\textbf v})=v_2+v_1v_2+\cdots.\]

\section{Another Matrix Approach}

For any given order $n$, the polynomials of degree $n$ are an invariant subspace for
the operator $Y$ up until the last step. We can formulate an alternative matrix
computation as follows. Let $\bar D$ and $\bar X$ denote the matrices of the operators of
differentiation and multiplication by $x$ respectively on polynomials of degree less than
or equal to $n$. The space is invariant under differentiation, and we cut off multiplication
by $x$ to be zero on $x^n$. We get
$$ \bar D_{ij}=i\,\delta_{i+1,j} \quad\hbox{and}\quad \bar X_{ij}=\delta_{i-1,j}$$
with the first row of $\bar X$ all zeros. We then compute the matrix $\bar X$ times
$W(\bar D)$, where $W(\bar D)$ is computed as a matrix polynomial by substituting in $W(z)$
up to order $n$. Then $Y$ has a matrix representation, $\bar Y=\bar X W(\bar D)$,
on the space and we iterate multiplying by $\bar Y$ acting on the unit vector $\bf e_1$. 
These give the coefficients of the polynomials $y_n$.

In several variables, one constructs matrices for $D_j$ and $X_i$ using Kronecker products
of $\bar D$ and $\bar X$ with the identity. For example, 
$$\bar D_j=I\otimes I\otimes\cdots\otimes\bar D\otimes I\cdots\otimes I$$ 
with $\bar D$ in the $j^{\rm th}$ spot. Similarly for $\bar X_i$.
Then one has explicit matrix representations for the
dual vector fields and the polynomials can be found accordingly.

This approach is explicit, but seems to much slower than using the built-in Ore\_algebra
package. 

\section{Worksheets}
\includegraphics[height=12.65cm,width=0.90\textwidth]{out1.epsi}

\includegraphics[totalheight=0.98\textheight,width=0.95\textwidth]{out2.epsi}

\includegraphics[scale=0.6]{out3.epsi}

\end{document}